\documentclass[a4paper]{amsart}
\usepackage{amssymb,amsmath,amsthm,amscd,epsfig}

\newtheorem{theorem}{Theorem}
\theoremstyle{plain}

\newtheorem{proposition}[theorem]{Proposition}

\title{On finite index subgroups of a universal group}
\author [Brumfiel,  Hilden,  Lozano,  Montesinos,
Ramirez, Short, Tejada,  Toro]{G. Brumfiel, H. Hilden, M.T.
Lozano*, J.M. Montesinos--Amilibia*, E. Ramirez--Losada, H. Short,
D. Tejada, D. Toro}
\address[G.~Brumfiel]{Stanford University, Stanford, Ca, USA}
\email[G.~Brumfiel]{brumfield@math.stanford.edu}
\address[H.~Hilden]{University of Hawaii, Honolulu, Hi, USA}
\email[H.~Hilden]{mike@math.hawaii.edu}
\address[M.T.~Lozano]{Universidad de Zaragoza, Zaragoza, Spain}
\email[M.T.~Lozano]{tlozano@unizar.es}
\address[J.M.~Montesinos]{Universidad Complutense, Madrid, Spain}
\email[J.M.~Montesinos]{montesin@mat.ucm.es}
\address[E.~Ramirez]{CIMAT, Mexico} \email{kikis@cimat.mx}
\address[H.~Short]{Universite de Provence, Marseille, France}
 \email{hamish.short@cmi.uni-mrs.fr}
 \address[D.~Tejada, M.~Toro]{Universidad Nacional de Colombia, Madellin, Colombia}
 \email{dtejada@unalmed.edu.co, mmtoro@unalmed.edu.co}
 \keywords{3-manifold, branched covering, universal link, universal group}
 \subjclass[2000]{57M12,57M25,57M50,57M60}
\thanks {$^*$This research
was supported by grants MTM2004088080 and MTM2006-00825}

\begin{document}
\begin{abstract}
The orbifold group of the Borromean rings with singular angle 90
degrees, $U$, is a universal group, because every closed oriented
3--manifold $M^{3}$ occurs as a quotient space $M^{3} = H^{3}/G$,
where $G$ is a finite index subgroup of $U$. Therefore, an
interesting, but quite difficult problem, is to classify the
finite index subgroups of the universal group $U$. One of the
purposes of this paper is to begin this classification. In
particular we analyze the classification  of the finite index
subgroups of $U$ that are generated by rotations.
\end{abstract}

\maketitle

\section{Introduction}

A finite covolume discrete group of isometries of hyperbolic 3--space, $%
H^{3} $, is said to be \emph{universal} if every closed oriented
3--manifold $M^{3}$ occurs as a quotient space $M^{3} = H^{3}/G$,
where $G$ is a finite index subgroup of the universal group. It
was originally shown in \cite{HLMW1987} that $U$, the orbifold
group of the Borromean rings with singular angle 90 degrees is
universal. (See \cite{BHLMSTT2006} for a simpler proof.)

Although there appear to be infinite families of universal groups, the group
$U$ is the only one so far known that is associated to a tessellation of $%
H^{3}$ by regular hyperbolic polyhedra in that there is a tessellation of $%
H^{3}$ by regular dodecahedra with dihedral angles $90^\circ$ any one of
which is a fundamental domain for $U$.

An interesting, important, but quite difficult problem, is to
classify the finite index subgroups of $U$. A theorem of Armstrong
\cite{A1968} shows that $\pi_1(M^{3}) \cong G/TOR(G)$ where
$TOR(G)$ is the subgroup of $G$ generated by rotations. In
particular $M^{3}$ is simply connected if and only if $G$ is
generated by rotations. One of the purposes of this paper is to
begin the classification of the finite index subgroups of $U$ that
are generated by rotations. Our main result is Theorem \ref{t7}.

\textbf{Theorem \ref{t7}} \textit{For any integer $n$ there is an
index $n$ subgroup of $U$ generated by rotations.}

In Theorem \ref{t8} we illustrate the essential differences
between the cases $n$ is odd and $n$ is even.

The organization of the paper is as follows: In Section 2 we define the
group $U$, a closely related Euclidean crystallographic group $\widehat U$,
and a homomorphism $\varphi : U \longrightarrow \widehat U$. In Section 3 we
show there are tessellations of $H^{3}$ by regular dodecahedra and $E^{3}$
by cubes and we exploit the homomorphism $\varphi : U \longrightarrow
\widehat U$ to define a branched covering space map $p : H^3 \longrightarrow
E^3$ that respects the two tesselations in the sense that the restriction of
$p$ to any one dodecahedron of the tesselation of $H^{3}$ is a homeomorphism
onto a cube of the tesselation of $E^{3}$. In Section 4 we prove the
\emph{rectangle theorem} and we use it to classify the finite index subgroups of $%
\widehat U$ that are generated by rotations. In the final section
we use this classification together with the homomorphism defined
in Section 2 to prove the main theorem of the paper, Theorem
\ref{t7}, and some existence theorems about finite index subgroups
of $U$ generated by rotations.

\section{Definitions of $U$, $\widehat U$ and the homomorphism $\protect%
\varphi : U \longrightarrow \widehat U$}

\begin{figure}[ht]
\epsfig{file=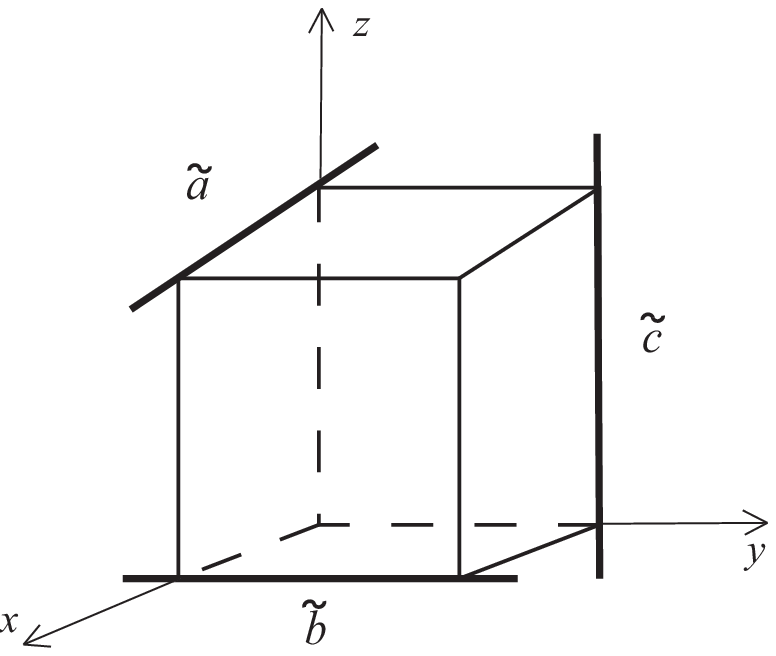,height=3.5cm}\caption{}\label{fig1}
\end{figure}

Let $C_{0}$ be the cube in $E^{3}$ with vertices $(\pm 1,\pm 1,\pm 1)$. We
obtain a tessellation of $E^{3}$ by applying compositions of even integer
translations in the $x$, $y$, and $z$ directions to $C_{0}$. In this paper
we do not consider any other tessellations of $E^{3}$ and we refer to this
tessellation as ``the'' tessellation of $E^{3}$. The intersection of $C_{0}$
with the positive octant, together with the lines $\widetilde{a}=(t,0,1)$, $%
\widetilde{b}=(1,t,0)$, and $\widetilde{c}=(0,1,t)$; $-\infty <t<\infty $,
is depicted in Figure \ref{fig1}.

The group $\widehat U$ is the Euclidean crystallographic group generated by
180 degree rotations $a$, $b$, and $c$ with axes $\widetilde a$, $\widetilde
b$, and $\widetilde c$, respectively. We see that $\widehat U$ preserves the
tessellation and contains the translations $t_x = b(cbc^{-1})$, $t_y =
a(cac^{-1})$, $t_z = a(bab^{-1})$, by distances of four, in the $x$, $y$,
and $z$ directions, respectively.

\begin{figure}[ht]
\epsfig{file=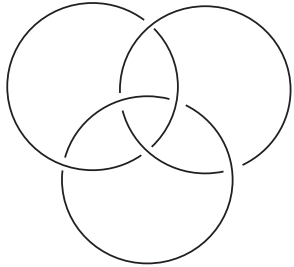,height=3.5cm}\caption{Borromean
rings.}\label{fig2}
\end{figure}

The cube $C_{0}$ is easily seen to be a fundamental domain for $\widehat{U}$%
, and the axes of rotation in $\widehat{U}$ divide each face of
each cube in the tessellation into two rectangles. The quotient
space $E^{3}/\widehat{U}$ is topologically $S^{3}$ as can be seen
by identifying faces of $C_{0}$ using $a$, $b$, $c$ and other
rotations. The group $\widehat{U}$ is the orbifold group of
$S^{3}$ as Euclidean orbifold with singular set the Borromean
rings $B$ and singular angle 180 degrees. This construction is due
to Thurston. For more details see (\cite{T1997},
\cite{BHLMSTT2006}). The Borromean rings are depicted in Figure
\ref{fig2}.

The induced map $p:E^{3}-\text{preimage\ }B\longrightarrow (E^{3}-\text{%
preimage\ }B)/\widehat{U}\approx S^{3}-B$ is a regular covering space map so
by the theory of covering spaces
\[
\widehat{U}\cong \pi _{1}(S^{3}-B)/p_{\ast }\pi _{1}(E^{3}-\text{preimage }%
B)\,.
\]
This gives rise to a presentation for $\widehat{U}$:
\[
\widehat{U}=\langle a,b,c| a\,b\overline{c}\overline{b}c=b\overline{c}%
\overline{b}c\,a,\,  b\,
c\overline{a}\overline{c}a=c\overline{a}\overline{c}a\,b,\, c\,
a\overline{b}\overline{a}b=a\overline{b}\overline{a}b\, c\,
,a^{2},\, b^{2},\,c^{2}\, \rangle \,.\tag{1}
\]
The presentation comes from the usual Wirtinger presentation of the group of
the Borromean rings with additional relations $a^{2}$, $b^{2}$, and $c^{2}$
arising from $p_{\ast }\pi _{1}$($E^{3}$--preimage $B$) which is normally
generated by squares of meridians about axes $\widehat{a}$, $\widehat{b}$,
and $\widehat{c}$.

There is a construction of $S^{3}$ as hyperbolic orbifold (also due to
Thurston) with singular set the Borromean rings analogous to the previous
construction. To describe it we shall work in the Klein model for $H^{3}$.

\begin{figure}[ht]
\epsfig{file=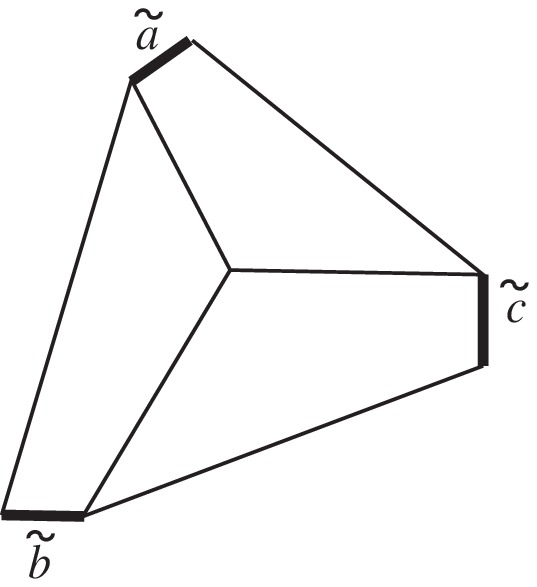,height=3.5cm}\caption{}\label{fig3}
\end{figure}

In the Klein model hyperbolic points are Euclidean points inside a ball of
radius $R$ centered at the origin in $E^{3}$ and hyperbolic lines and planes
are the intersections of Euclidean lines and planes with the interior of the
ball of radius $R$. Let $D_{0}$ be a regular Euclidean dodecahedron that is
symmetric with respect to reflection in the $xy$, $yz$, and $xz$ planes. The
intersection of $D_{0}$ with the positive octant is depicted in Figure \ref
{fig3}.

If $R$ is chosen correctly, (Details are in \cite{HLM2001}), then $%
D_0 $ can be considered as a regular hyperbolic dodecahedron with 90 degree
dihedral angles. Each pentagonal face contains one edge that lies in either
the $xy$, $xz$, or $yz$ plane. Reflection in this plane, restricted to the
pentagon, defines an identification in pairs on the pentagonal faces of $D_0$%
. As in the construction with the cube $C_0$, the resulting topological
space is $S^{3}$. A hyperbolic orbifold structure is thus induced on $S^{3}$
with singular set the Borromean rings, $B$, and singular angle 90 degrees.
The Borromean rings are the image, after identification of the pentagonal
edges that lie in the $xy$, $xz$, and $yz$ planes.

There is a 4--fold regular branched cyclic covering
$q_{1}:X^{3}\longrightarrow S^{3}$ with branch set the Borromean
rings induced by the natural group homomorphisms
\[\pi _{1}(S^{3}-B)\longrightarrow H_{1}(S^{3}-B;Z)\cong
Z\oplus Z\oplus Z\longrightarrow Z\text{ mod} 4.\]
 The hyperbolic
orbifold structure on $S^{3}$ with singular set the Borromean
rings pulls back to a hyperbolic manifold (not orbifold) structure
on $X^{3}$ as meridians are sent to 1 in the above homomorphism.

The hyperbolic manifold $X^{3}$ has a \emph{tessellation}
consisting of four dodecahedra each of which is sent
homeomorphically to $D_{0}$ by the map $p$. The universal covering
space map $q_{2}:H^{3}\longrightarrow X^{3}$ is used to pull back
the tessellation of $X^{3}$ by dodecahedra to a tessellation of
$H^{3}$ by dodecahedra. The composition of covering space maps
$q_{1}\circ
q_{2}:H^{3}\longrightarrow S^{3}$ is a regular branched covering space map $%
H^{3}\longrightarrow S^{3}$ induced by the group of hyperbolic isometries $U$%
. That is to say there is a quotient branched covering map $%
H^{3}\longrightarrow H^{3}/U\approx S^{3}$ and an associated unbranched
covering space map $p:H^{3}-\text{axes of rotation }=H^{3}-\text{preimage }%
B\longrightarrow H^{3}-\text{preimage }B/U\approx S^{3}-B$. As in the
Euclidean case this covering space map gives rise to a presentation for $U$
via covering space theory:
\[
U=\langle a,b,c| a\,b\overline{c}\overline{b}c=b\overline{c}%
\overline{b}c\,a,\,  b\,
c\overline{a}\overline{c}a=c\overline{a}\overline{c}a\,b,\, c\,
a\overline{b}\overline{a}b=a\overline{b}\overline{a}b\, c\,
,a^{4},\, b^{4},\,c^{4}\, \rangle \tag{2}
\]

As before the presentation comes from the usual Wirtinger
presentation of the group of the Borromean rings with additional
relations $a^{4}$, $b^{4}$, $c^{4}$ arising from $p_{\ast }\pi
_{1}(H^{3}-\text{preimage} B)$  which is normally generated by
fourth powers of meridians about the axes $ \widetilde{a}$,
$\widetilde{b}$ and $\widetilde{c}$.

Examining the presentations for $U$ and $\widehat U$ we see that they are
the same except for the relations $a^4$, $b^4$, and $c^4$ in $U$ and $a^2$, $%
b^2$, and $c^2$ in $\widehat U$. Nonetheless the map $a \to a$, $%
b \to b$, and $c \to c$, mapping generators of $U$ to generators
of $\widehat U$, defines a homomorphism $\varphi : U
\longrightarrow \widehat U$ and an exact sequence.
\[
1 \longrightarrow K \longrightarrow U
\overset{\varphi}{\longrightarrow} \widehat U \longrightarrow 1\,.
\tag{3}
\]
In this exact sequence $K$ is defined to be the kernel of
homomorphism $ \varphi$.

We say that a group of isometries of $H^{3}$ or $E^{3}$ is
\emph{associated} to a tessellation of $H^{3}$ or $E^{3}$ by
regular compact polyhedra if there is a tessellation of $H^{3}$ or
$E^{3}$ by regular compact polyhedra any one
of which is a fundamental domain for the group. Thus the groups $U$ and $%
\widehat{U}$ are associated to the tessellations of $H^{3}$ and $E^{3}$ by
regular dodecahedra and cubes, respectively. This is not a common
occurrence. For example, of the regular polyhedra only cubes can tessellate $%
E^{3}$. In the table below, we have listed the cosines of the dihedral
angles of the Euclidean regular polyhedra and also the dihedral angles of
the hyperbolic regular polyhedra with vertices on the sphere at infinity.
Tetrahedra, octahedra, dodecahedra and icosahedra cannot tessellate $E^{3}$
because their dihedral angles are not submultiples of $360^{\circ }$ so they
don't ``fit around an edge''.
\medskip

\begin{tabular}{|c|c|c|}\hline
Polyhedral Type & Euclidean  dihedral angle& Hyperbolic
dihedral angle  \\
& & vertices at $\infty $ \\
 \hline Tetrahedron & $ArcCos[1/3]\approx 70.5288^{\circ }$ &
$60^{\circ }$
 \\ \hline
Cube &ArcCos[0]= $90^{\circ }$& $60^{\circ }$
\\ \hline
Octahedron & $ArcCos[-1/3]\approx 109.471^{\circ }$  & $90^{\circ
}$
\\ \hline
Dodecahedron & $ArcCos[-1/\sqrt{5}]\approx 116.565^{\circ }$&
$60^{\circ }$
\\ \hline
Icosahedron & $ArcCos[-\sqrt{5}/3]\approx 138.19^{\circ }$ &
$108^{\circ }$ \\ \hline
\end{tabular}
\medskip

There are five regular Euclidean polyhedra but the corresponding hyperbolic
polyhedra occur in one parameter families. One can construct the family of
hyperbolic cubes, for example, by starting with $C_0$, the cube with
vertices $(\pm1, \pm1, \pm1)$, in the Klein model with the sphere at
infinity having Euclidean radius $R = \sqrt3$ and let $R$ increase from $%
\sqrt3$ to $\infty$. There is an isometry from the Klein model using the
Euclidean ball of radius $R$ to the Poincar\'e model using the same
Euclidean ball (as Thurston has explained), that is the identity on the
sphere at infinity. Since the Poincar\'e model is conformal and Poincar\'e
hyperbolic planes are Euclidean spheres perpendicular to the sphere at
infinity, the dihedral angle between two Poincar\'e planes is the same as
the Euclidean angle between the two circles in which the Poincar\'e planes
intersect the sphere at infinity. Thus the dihedral angle between two Klein
planes is the same as the angle between the two circles in which they
intersect the sphere at infinity. As $R$ increases, in the case of the cube,
for example, from $\sqrt3$ to infinity, the dihedral angle increases from $%
60^\circ$ to $90^\circ$. There exists a compact hyperbolic cube with
dihedral angle $\theta$ if and only if $0 < \cos\theta < 1/2$. Thus, if it
is possible to tessellate $H^{3}$ with compact hyperbolic cubes they must
have dihedral angle 72 degrees as that is the only submultiple of $360^\circ$
in the range of possible dihedral angles. A glance at the table (4)
indicates that it is impossible to tessellate $H^{3}$ with compact regular
octahedra or tetrahedra and if it is possible to tessellate $H^{3}$ with
icosahedra the dihedral angle must be 120 degrees. In the dodecahedral case
we have shown that there is a tessellation of $H^{3}$ by regular compact
hyperbolic dodecahedra with dihedral angle 90 degrees. If there were a
different tessellation by compact regular dodecahedra the dihedral angle
would have to be $72^\circ$.

All the above is part of standard 3--dimensional hyperbolic geometry and we
explain it mainly so as to highlight the singular nature of the groups $U$
and $\widehat U$ and the tessellations with which they are associated and as
background for the following conjecture.

\textbf{Conjecture} \textit{The group $U$ is the only universal
group associated to a tessellation of $H^{3}$ by regular
hyperbolic polyhedra.}

In the next section we study the groups $U$ and $\widehat U$ and
the tessellations to which they are associated to produce a
branched covering of $E^{3}$ by $H^{3}$.

\section{$H^3$ as a branched covering of $E^3$}

Let $D_{0}$ and $C_{0}$ be the regular dodecahedron and cube in the Klein
model for $H^{3}$ and in $E^{3}$ respectively, as defined in the previous
section. We know that $D_{0}$ is a fundamental domain for the group $U$ and
is also an element of the tessellation of $H^{3}$ by regular dodecahedra.
For any other dodecahedron $D$ in the tessellation there is a unique element
$u$ of $U$ such that $u(D_{0})=D$. Analogously, $C_{0}$ is a fundamental
domain for the group $\widehat{U}$ and is part of the tessellation of $E^{3}$
by cubes. For any other cube $C$ in the tessellation there is a unique
element $\widehat{u}$ of $\widehat{U}$ such that $\widehat{u}(C_{0})=C$.

Let $\alpha _{0}:D_{0}\longrightarrow C_{0}$ be a homeomorphism that is as
nice as possible. Thus $\alpha _{0}$ should commute with reflections in the $%
xy$, $xz$, and $yz$ planes and also with the 3--fold rotations about the
axes $\{(t,t,t)\}$ in the Klein model for $H^{3}$ and in $E^{3}$. The cube $%
C_{0}$ becomes a dodecahedron when each of its faces is split in half by an
axis of rotation of $\widehat{U}$. Then $\alpha _{0}$, viewed as a map
between dodecahedra takes vertices, edges, and faces to vertices, edges, and
faces, respectively.

Now we define a map $p:H^{3}\longrightarrow E^{3}$. Let $p=\alpha _{0}$ on $%
D_{0}$. Any other point $A$ in $H^{3}$ belongs to a dodecahedron $D$ of the
tessellation. There is a unique $u\in U$ such that $u(D_{0})=D$. Let $%
\widehat{u}=\varphi (u)$ where $\varphi :U\longrightarrow \widehat{U}$ is
the homomorphism defined in the previous section. Define the map $p$ by $%
p(A)=\widehat{u}\circ \alpha _{0}\circ u^{-1}(A)$. The map $p$ is well
defined for points in the interior of dodecahedra in the tessellation but we
must show that $p$ is well defined for the other points. Let $A$ belong to
the interior of a pentagonal face $P$ belonging to each of two adjacent
dodecahedra $D_{1}$ and $D_{2}$.

Then there are unique elements $u_{1}$ and $u_{2}$ of $U$ such that $%
u_{1}(D_{0})=D_{1}$ and $u_{2}(D_{0})=D_{2}$. Then $u_{1}^{-1}(D_{2})$ is a
dodecahedron, call it $\widehat{D}$, that intersects $D_{0}$ exactly in a
pentagonal face $P_{0}$. The pentagonal face $P_{0}$ of $D_{0}$ intersects
exactly one of the six axes of rotation, call it $ax$, that intersect $D_{0}$
and this axis lies in the $xy$, $xz$, or $yz$ plane of the Klein model.
There is a $90^{\circ }$ rotation about $ax$, call if $rot$, that sends $%
D_{0}$ to $\widehat{D}$. Thus $u_{1}\circ rot(D_{0})=D_{2}$ which implies $%
u_{1}\circ rot=u_{2}$, which further implies $\widehat{u}_{1}\circ \widehat{%
rot}=\widehat{u}_{2}$ in group $\widehat{U}$. Then $\widehat{u}_{2}\circ
\alpha _{0}\circ u_{2}^{-1}=\widehat{u}_{2}\circ \alpha _{0}\circ
rot^{-1}\circ u_{1}^{-1}=\widehat{u}_{1}\circ \widehat{rot}\circ \alpha
_{0}\circ rot^{-1}\circ u_{1}$ so that to show that the map $p$ is well
defined on the interior of pentagon $P$ it suffices to show that $\widehat{%
rot}\circ \alpha _{0}\circ rot^{-1}=\alpha _{0}$ when restricted to
pentagonal face $P_{0}$.

The homomorphism $\varphi :U\longrightarrow \widehat{U}$ takes $a$, $b$, and
$c$ to $\widehat{a}$, $\widehat{b}$, $\widehat{c}$, respectively where $a$, $%
b$, and $c$ are $90^{\circ }$ rotations about axes $\widetilde{a}$, $%
\widetilde{b}$, and $\widetilde{c}$, respectively of Figure \ref{fig3} and $\widehat{a%
}$, $\widehat{b}$, $\widehat{c}$ are $180^{\circ }$ rotations about axes $%
\widetilde{a}$, $\widetilde{b}$, and $\widetilde{c}$, respectively
of Figure \ref{fig1}. The rotation $rot$ is one of $a$, $b$, $c$,
$a^{-1}$, $b^{-1}$, $c^{-1}$,
$bab^{-1}$, $cbc^{-1}$, $aca^{-1}$, $ba^{-1}b^{-1}$, $cb^{-1}c^{-1}$, $%
ac^{-1}a^{-1}$. The rotation $rot$, when restricted to pentagon $P_{0}$
equals reflection in the $xy$, $yz$, or $xz$ plane depending on which plane
axis $rot$ lies in. Similarly, the rotation $\widehat{rot}$ is one of $%
\widehat{a}$, $\widehat{b}$, $\widehat{c}$, $\widehat{b}\widehat{a}\widehat{b%
}^{-1}$, $\widehat{c}\widehat{b}\widehat{c}^{-1}$, $\widehat{a}\widehat{c}%
\widehat{a}^{-1}$ and the rotation $\widehat{rot}$ when restricted to the
half square that is the image of $P_{0}$ under $\alpha $ equals reflection
in the $xy$, $xz$, or $yz$ plane depending on which plane axis $\widehat{rot}
$ lies in. But $\alpha _{0}$ commutes with reflections in the $xy$, $xz$, or
$yz$ planes so that $\widehat{rot}\circ \alpha _{0}\circ rot^{-1}=\alpha
_{0} $ and the map $p$ is well defined on the interiors of dodecahedra in
the tessellation and on the interiors of their pentagonal faces. That $p$ is
also well defined on edges and vertices of the tessellating dodecahedra now
follows by a continuity argument.

We summarize all this in a theorem.

\begin{theorem}
\label{t1} There exists a tessellation of $H^{3}$ by regular
hyperbolic dodecahedra with $90^{\circ}$ dihedral angle and a
tessellation of $E^{3}$ by cubes and a map $p:H^{3}\longrightarrow
E^{3}$ such that the following holds.

1. Any dodecahedron in the tessellation of $H^{3}$ is a fundamental domain
for the universal group $U$.

2. Any cube in the tessellation of $E^{3}$ is a fundamental domain for the
Euclidean crystallographic group $\widehat{U}$.

3. The axes of rotation in $\widehat{U}$ divide each face of each cube in
the tessellation of $E^{3}$ into two rectangles so that the cube may be
viewed as a dodecahedron.

4. The restriction of $p$ to any one dodecahedron is a homeomorphism of that
dodecahedron onto a cube in the tessellation of $E^{3}$. When the cube is
viewed as a dodecahedron as in 3 above, the map $p$ sends vertices, edges,
and faces to vertices edges and faces respectively. The map $p$ also sends
axes of rotation for $U$ homeomorphically, even isometrically, to axes of
rotation for $\widehat{U}$.

5. The map $p$ is a branched covering space map with all branching of order
two.
\end{theorem}

In effect, parts 1 through 4 of the theorem have already been
proven in the remarks preceding the statement of the theorem. To
see that 5 is true, it is only necessary to examine $p$ near an
axis of rotation for $U$. The branching is of order two because
four dodecahedra \emph{fit} around every axis of rotation in $U$
while only two cubes \emph{fit }around an axis of rotation of
$\widehat{U}$.

It is clear from the definition of the map $p$ when restricted to a
dodecahedron, $p=\widehat{u}\circ \alpha _{0}\circ u^{-1}$, that the group
of covering transformations is the kernel of the homomorphism $\varphi
:U\longrightarrow \widehat{U}$. On the other hand $p$ when restricted to ($%
H^{3}$--axes of rotation for $U$) is an unbranched covering of ($E^{3}$ --
axes of rotation for $\widehat{U}$) so that $K=\ker \varphi
:U\longrightarrow \widehat{U}$ is isomorphic to $\pi _{1}(E^{3}$ -- axes of
rotation for $\widehat{U}$) modulo $p_{\ast }\pi _{1}$($H^{3}$ -- axes of
rotation for $U$), by standard covering space theory.

As $\pi _{1}(E^{3}$ -- axes of rotation for $\widehat{U}$) is a free group
generated by meridians, one meridian for each axis of rotation, and $\pi
_{1}(H^{3}$ -- axes of rotation for $U$) is also generated by meridians it
follows that $p_{\ast }\pi _{1}(H^{3}$ -- axes of rotation) is normally
generated by squares of meridians, one for each axis of rotation in $%
\widehat{U}$. We also summarize all this in a theorem.

\begin{theorem}
\label{t2} The group of covering transformations for the branched
covering \linebreak $p:H^{3}\longrightarrow E^{3}$ is isomorphic
to the group $K$ that is the kernel of $\varphi :U\longrightarrow
\widehat{U}$. The group $K$ is naturally isomorphic to a countable
free product of $Z$ mod 2's, one generator for each axis of
rotation in $\widehat{U}$. In particular the group $K$ is
generated by 180 degree rotations.
\end{theorem}

As before, the proof of the theorem is in effect given by the remarks
immediately prior to the statement of the theorem.

Theorems \ref{t1} and \ref{t2} enable us to ``label'' each axis of
rotation in $U$ with an algebraic integer in the field
$Q(\sqrt{-3}\,)$. Note that each axis of rotation for
$\widehat{U}$ is a line of parametric equation ($t$, even, odd) or
(odd, $t$, even) or
(even, odd, $t$), $-\infty <t<\infty $. Any such axis intersects the plane $%
x+y+z=0$ in a point (odd, odd, even) or (even, odd, odd) or (odd,
even, odd) as zero is even. One can verify that the intersection
of the tessellation by cubes of $E^{3}$ with the plane $x+y+z=0$
induces a tessellation of the plane $\pi :\,x+y+z=0$ by (regular)
hexagons and (equilateral) triangles and that
cube $C_{0}$ intersects the plane $x+y+z=0$ in a hexagon with vertices $%
\{(\pm 1,\mp 1,0),(\pm 1,0,\mp 1),(0,\pm 1,\mp 1)\}$. Using a
similarity of the plane $x+y+z=0$ with center the origin and
expansion ratio $1/\sqrt{2}$ we can recoordinatize the plane
$x+y+z=0$ by the complex numbers $\mathbb{C}$ so that the six
vertices of this hexagon have coordinates equal to the six roots
of unity in $\Bbb{C}$. Then every axis of rotation of
$\widehat{U}$ intersects the plane $x+y+z=0$ in a point whose
coordinate is an algebraic integer in the field $Q(\sqrt{-3}\,)$.
We label each axis $d$ of rotation of $U$ with the coordinate of
$p(d)\cap \pi$. Again we summarize these results in a theorem.

\begin{theorem}
\label{t3} In the branched covering $p:H^{3}\longrightarrow E^{3}$ each axis
of rotation for $U$ is labelled by an algebraic integer of the field $Q(%
\sqrt{-3}\,)$. The group of covering transformations $K$ preserves
labelling. For any two axes of rotation $a$ and $b$ of $U$ with the same
label, there is an element $k$ of $K$ such that $k(a)=b$.
\end{theorem}

In the next section we classify the subgroups of finite index in $\widehat U$
that are generated by rotations.

\section{Finite index subgroups of $\widehat U$ generated by rotations}

The group $\widehat{U}$ is the crystallographic group
$I2_{1}2_{1}2_{1}$, number 24 of the International Tables of
Crystallography \cite{tabcrys} . In this section we describe two
families of subgroups of $\widehat{U}$ (defined in Section 2)
generated by rotations. And we show that \textit{any} finite index
subgroup of $\widehat{U}$ generated by rotations is equivalent (in
a sense we make precise) to exactly one member of one of the two
families.

The axes of rotation of $\widehat{U}$ have parametric equations of
form ($t$, even, odd), (odd, $t$, even) or (even, odd, $t$);
$-\infty <t<\infty $ according as to whether they are parallel to
the $x$, $y$, or $z$ axes. The distance between axes lying in a
plane parallel to the $xy$, $xz$, or $yz$ planes is an even
integer.

\begin{figure}[ht]
\epsfig{file=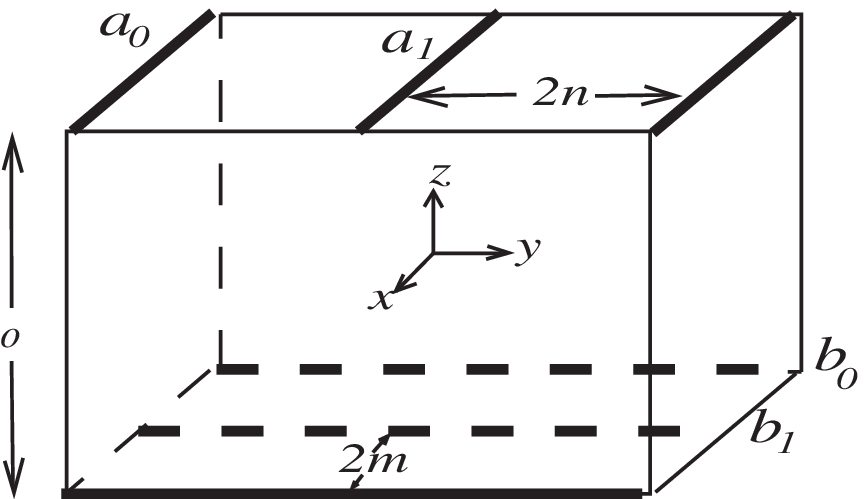,height=3.5cm}\caption{$Box(G(m,n,o)$.}\label{fig4}
\end{figure}

Let $(m,n,o)$ be a triple of positive integers where $o$ is odd and $m$ and $%
n$ are arbitrary. Let $Box(G(m,n,o))$  be the rectangular
parallelepiped defined by the following conditions.

a. The front and back faces of $Box(G(m,n,o))$ lie in the planes $%
x=2m+1 $ and $x=-2m+1$, respectively.

b. The right and left faces of $Box(G(m,n,o))$ lie in the planes
$y=2n$ and $y=-2n$, respectively.

c. The top and bottom of $Box(G(m,n,o))$ lie in the planes $z=o$ and $%
z=0$, respectively. $Box(G(m,n,o))$ together with certain axes of
rotation is pictured in Figure \ref{fig4}.

Axes $a_{0}$, $a_{1}$, $b_{0}$ and $b_{1}$ have parametric equations $%
(t,-2n,o)$, $(t,0,o)$, $(-2m+1,t,0)$ and $(1,t,0)$, respectively. Then $%
\widehat{G}(m,n,o)$ is defined to be the subgroup of $\widehat{U}$ generated
by $A_{0}$, $A_{1}$, $B_{0}$, and $B_{1}$, the rotations in the axes $a_{0}$%
, $a_{1}$, $b_{0}$, and $b_{1}$, respectively.

Observe that $T_{y}=A_{1}A_{0}$, $T_{x}=B_{1}B_{0}$, and $%
T_{z}=(A_{0}B_{1})^{2}$ are translations by $4m$, $4n$, and $4o$ in the $x$ $%
y$, and $z$ directions, respectively. Another generating set of $\widehat{G}%
(m,n,o)$ is $A_{1}$, $B_{1}$, $T_{x}$, and $T_{y}$. Conjugating a
translation $T_{x}$, $T_{y}$, or $T_{z}$ by a rotation $A_{1}$,
$B_{1}$ either results in the translation itself or its inverse,
so there are commutation relations such as
$B_{1}T_{x}=T_{x}^{-1}B_{1}$. Thus any element of
$\widehat{G}(m,n,o)$ has form $T$, $A_{1}T$, $B_{1}T$ or
$A_{1}B_{1}T$
where $T$ is a translation that is some product of $T_{x}$, $T_{y}$, and $%
T_{z}$. With these observations we can see that $Box(G(m,n,o))$ is
a fundamental domain for the group $\widehat{G}(m,n,o)$. The
volume of $Box(G(m,n,o))$ equals $4m\times 4n\times o$ and the
volume of cube $C_{0}$, which is a fundamental domain for
$\widehat{U}$ equals 8. Thus dividing one
by the other, the index of $\widehat{G}(m,n,o)$ in $\widehat{U}$ equals $%
2mno $, an even integer. The group $\widehat{G}(m,n,o)$ is the
crystallographic  group $P222_{1}$, number 17 in the International
Tables of Crystallography \cite{tabcrys}.

Let $(p,q,r)$ be a triple of odd positive integers such that $p\leqq q$ and $%
p\leqq r$ and if the three integers are not all different then
$p\leqq q\leqq r$. Let $Box(H(p,q,r))$ be the rectangular
parallelepiped defined by the following conditions.

\begin{figure}[ht]
\epsfig{file=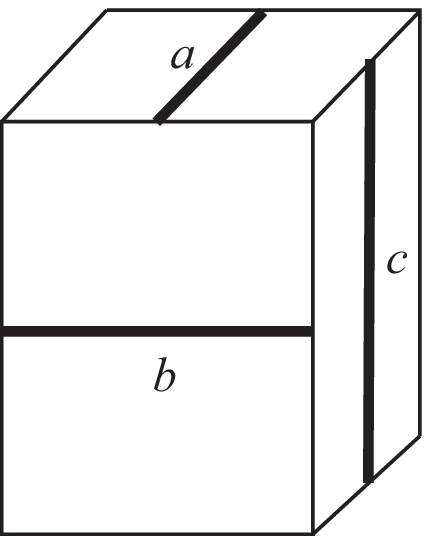,height=3.5cm}\caption{$Box(H(p,q,r))$.}\label{fig5}
\end{figure}

The front and back, left and right, top and bottom faces of
$Box(H(p,q,r))$ lie in the planes $x=p$, $x=-p$; $y=q$, $y=-q$;
$z=r$, $z=-r$, respectively. $Box(H(p,q,r))$ is pictured in Figure
\ref{fig5} along with axes of rotation $a=(t,0,r)$, $b=(p,t,0)$,
and $c=(0,q,t)$.

The group $\widehat{H}(p,q,r)$ is defined to be the subgroup of
$\widehat{U}$ generated by rotations $A$, $B$, and $C$ in axes
$a$, $b$, and $c$, respectively. Observe that $T_{x}=BCBC$,
$T_{y}=CACA$, and $T_{z}=ABAB$ are translations by $2p$, $2q$, and
$2r$ in the $x$, $y$, and $z$ directions, respectively. Also note
that conjugating $T_{x}$, $T_{y}$, or $T_{z}$ by ($A$
or $B$ or $C$) results in $T_{x}$ or $T_{x}^{-1}$, $T_{y}$ or $T_{y}^{-1}$, $%
T_{z}$ or $T_{z}^{-1}$, respectively. These observations imply
that any element of group $\widehat{H}(p,q,r)$ equals exactly one
of $T$, $AT$, $BT$ or $CT$ where $T$ is a product of $T_{x}$,
$T_{y}$ and $T_{z}$. As before, we can see that $Box(H(p,q,r))$ is
a fundamental domain for group $\widehat{H}(p,q,r)$. The group
 $\widehat{H}(p,q,r)$ is again the crystallographic  group $I2_{1}2_{1}2_{1}$,
number 24 in \cite{tabcrys}.

The volume of $Box(H(p,q,r))$ equals $8pqr$ and volume $C_{0}=8$
so,
reasoning as before, the index of $\widehat{H}(p,q,r)$ in $\widehat{U}$ is $%
pqr$ which is an odd integer.

We wish to define an equivalence relation on infinite index subgroups of $%
\widehat{U}$. Let $D$ be the $120^{\circ }$ rotation about the axis $(t,t,t)$%
; $-\infty <t<\infty $, which is a main diagonal of cube $C_{0}$ and let $%
\widehat{S}$ be the group generated by $D$ and $\widehat{U}$. As $D$ has
order three and normalizes $\widehat{U}$ we see that $[\widehat{S}:\widehat{U%
}]=3$. We define two subgroups of $\widehat{U}$ to be equivalent if they are
conjugate as subgroups of $\widehat{S}$. This equivalence relation leads to
the least messy classification of the finite index subgroups of $\widehat{U}$
generated by rotations. We observe that rotation $D$ cyclically permutes the
$x$, $y$, and $z$ axes but that there is no element of $\widehat{S}$ that
fixes one of these three axes while transposing the other two.

The triple ``distance between adjacent axes'' parallel to the $x$, $y$, $z$
axes, respectively, defines an invariant on the groups $\widehat{G}(m,n,o)$
and $\widehat{H}(p,q,r)$. Thus triple $(\widehat{G}(m,n,o))=(2n,2m,\text{none%
})$ and triple $(\widehat{H}(p,q,r))=(2r,2p,2q)$.

Conjugating a $\widehat{G}$ or an $\widehat{H}$ by an element of $\widehat{S}
$ at most changes a triple by cyclically permuting it. Thus the fact that $%
\widehat{G}$ contains no axes of rotation parallel to the $z$--axis implies
that if $\widehat{G}(m,n,o)\sim \widehat{G}(\widetilde{m},\widetilde{n},%
\widetilde{o})$ then $(m,n,o)=(\widetilde{m},\widetilde{n},\widetilde{o})$
and the conditions $p\leqq q$ and $p\leqq r$, etc., imply that if $\widehat{H%
}(p,q,r)=\widehat{H}(\widetilde{p},\widetilde{q}\ \widetilde{r})$ then $%
(p,q,r)=(\widetilde{p},\widetilde{q},\widetilde{r})$. Also as the index of a
$\widehat{G}$ in $\widehat{U}$ is even and the index of an $\widehat{H}$ in $%
\widehat{U}$ is odd no $\widetilde{G}$ can be equivalent to an $\widehat{H}$%
. The rest of the classification consists of showing that any finite index
subgroup of $\widehat{U}$ generated by rotations is either equivalent to an $%
\widehat{H}$ or a $\widehat{G}$.

Suppose that $\widehat{G}$ is a finite index subgroup of $\widehat{U}$ that
is generated by rotations. If $\widehat{G}$ contained only rotations
parallel to one of the three axes, it would leave planes perpendicular to
this axis invariant and thus have infinite index in $\widehat{U}$. So $%
\widehat{G}$ either contains rotations about axes parallel to two of the
three axes $x$, $y$, and $z$ or it contains rotations about axes parallel to
all three. In the former case, we can assume $\widehat{G}$ contains
rotations with axes parallel to the $x$ and $y$ axes but doesn't contain
rotations with axes parallel to the $z$--axis by conjugating by an element
of $\widehat{S}$ if need be. In either case let $\mathcal{P}$ be a plane
parallel to the $yz$ plane in which an axis of $\widehat{G}$ parallel to the
$y$--axis lies. The set of axes of rotation of $\widehat{G}$ parallel to the
$x$--axis intersects $\mathcal{P}$ in a set of points we call axis points.

\begin{proposition}
\label{p4} (The rectangle theorem) There is a tessellation of $\mathcal{P}$
by congruent rectangles with sides parallel to the $y$ and $z$ axes such
that the set of axis points equals the set of vertices of the rectangles.
Each rectangle is divided in half by an axis of rotation for $\widehat{G}$
parallel to the $y$--axis.
\end{proposition}

The proof of Proposition 4 rests on three facts.

1. If $A$ is a rotation in $\widehat{G}$ with axis $\ell $ and $S\in
\widehat{G}$ then $SAS^{-1}$ is a rotation in $\widehat{G}$ with axis $%
S(\ell )$. In particular if $X$ is an axis point and $S(\mathcal{P})=%
\mathcal{P}$, then $S(X)$ is an axis point.

2. If $A$ is a rotation in $\widehat{G}$ with axis $\ell $ and $T$ is a
translation in $\widehat{G}$ such that $T(\mathcal{P})=\mathcal{P}$ and $%
\ell \cap \mathcal{P}=X$ then $TA$ is also a rotation in $\widehat{G}$ and
axis $(TA)\cap \mathcal{P}$ is the midpoint of the line segment $XT(X)$.

3. Group $\widehat{G}$ contains translations in the $x$, $y$, and $z$
directions. (Because $\widehat{U}$ does and $[\widehat{U}:\widehat{G}\
]<\infty $.)

\begin{figure}[ht]
\epsfig{file=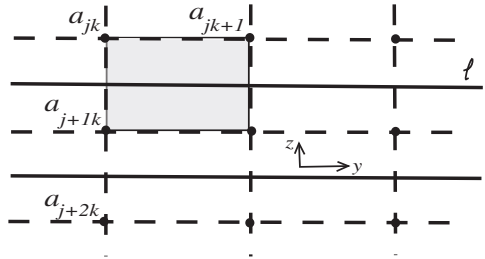,height=3.5cm}\caption{The plane
$\mathcal{P}$. }\label{fig6}
\end{figure}
\begin{proof}[Proof of Proposition 4]
Let $T_{y}$ and $T_{z}$ be translations by minimal distance in the
$y$ and $z $ directions respectively, belonging to $\widehat{G}$.
(Refer to Figure \ref{fig6}.) Let $a_{00}$ be an axis point. Then
by 1 and 2 above, $a_{20}=T_{y}(a_{00})$
and $a_{10}=$ midpoint $a_{00}a_{20}$ are axis points as are $%
a_{02}=T_{z}(a_{00})$, $a_{01}=$ midpoint $a_{00}a_{02}$ and $a_{11}=$
midpoint $a_{01}T_{y}(a_{01})$. The set of vertices of the tessellation by
rectangles referred to in Proposition 4 equals $\{T_{y}^{i}T_{z}^{j}a_{k\ell
}\mid i,j\in Z\ k,\ell \in \{0,1\}\}$.

Suppose $\ell $ is the axis of rotation of $B$ and $\ell $ lies in plane $%
\mathcal{P}$, is parallel to the $y$--axis and intersects the rectangle $%
R=\{a_{jk},a_{j+1k},a_{jk+1},a_{j+1k+1}\}$, where axis point $a_{jk}$
corresponds to rotation $A_{jk}$, etc. Then $\ell $ cannot contain the
vertices of $R$ as axes of rotation of distinct elements of $\widehat{U}$
don't intersect. And $\ell $ must divide $R$ exactly in half for if $\ell $
lay closer to $a_{jk}$ than to $a_{j+1k}$ the element $A_{jk}(BA_{jk}B^{-1})$
of $\widehat{G}$ would be a translation in the $y$--direction by a distance
less than $a_{jk}a_{jk+2}$ contradicting the minimality in the choice of $%
T_{y}$. The set of translates of the axes $\ell $ and $A_{jk}(\ell )$ divide
every rectangle of the tessellation in half. We must show there are no axis
points in $\mathcal{P}$ not of the form $a_{jk}$. Suppose $x$ was such a
point corresponding to rotation $X$ and lying in rectangle $%
R=\{a_{jk},a_{jk+1},a_{j+1k},a_{j+1k+1}\}$. Then $x$ cannot lie on the sides
of the rectangle. (For example, if $x$ lay on $a_{jk}a_{jk+1}$, $%
XA_{yk}X^{-1}A_{jk}$ would be a translation in the $y$ direction by less
than length $a_{jk}a_{jk+2}$ contradicting the minimality in the choice of $%
T_{y}$.) And $x$ cannot lie on $\ell $. As $x$ belongs to the interior of
the rectangle and not on $\ell $, $X(BXB^{-1})$ is a translation in the $y$%
--direction by a distance less than $a_{jk}a_{jk+2}$ which is impossible.
\end{proof}

The next problem is to construct a fundamental domain for $\widehat{G}$.
With this in mind select a plane $\mathcal{P}$ parallel to the $yz$ plane
containing an axis $\ell $ in $\widehat{G}$ that is parallel to the $y$%
--axis. Recall that axes in $\widehat{G}$ parallel to the $x$, $y$, or $z$
axis have parametric equations ($t$, even, odd), (odd, $t$, even) or (even,
odd, $t$) respectively. Thus plane $P$ has equation $x=O$ where $O$ is odd.
Define the rectangle $R_{1}$ in $\mathcal{P}$, as pictured in Figure 7,
bounded on one side by $\ell $ with parametric equation $(O,t,e_{1})$ with $%
e_{1}$ even and having the opposite two vertices be axis points for $%
\mathcal{P}$ with coordinates $(O,E,o_{1})$ and $(O,E+4n,o_{1})$ with $o_{1}$
odd.

\begin{figure}[ht]
\epsfig{file=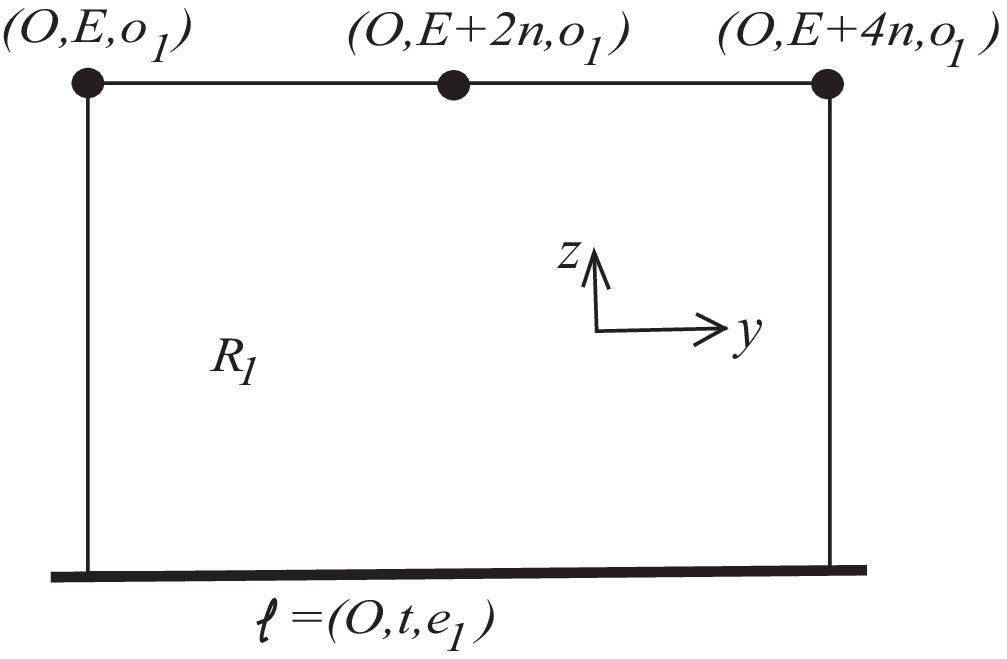,height=3.5cm}\caption{}\label{fig7}
\end{figure}
\begin{figure}[ht]
\epsfig{file=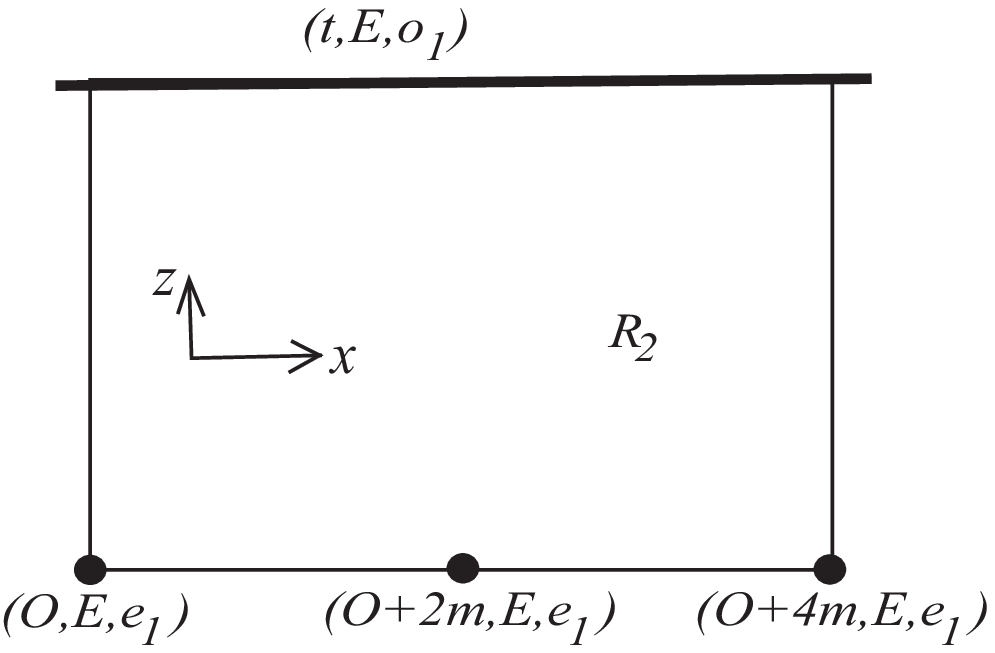,height=3.5cm}\caption{}\label{fig8}
\end{figure}

There is a rectangle theorem analogous to Proposition \ref{p4} but
with $x$ substituted for $y$. Let $\mathcal{Q}$ be the plane $y=E$
which contains the $x$ axis from $\widehat{G}$ with equation
$(t,E,o_{1})$. Then $\mathcal{Q}$ also is tesselated by rectangles
and we define $R_{2}$ to be the rectangle pictured in Figure
\ref{fig8}. Like $R_{1}$, the rectangle $R_{2}$ is not part of the
tessellation but is formed by gluing two half--rectangles from the
tessellation. $R_{2}$ is bounded on one side by axis $(t,E,o_{1})$
and the two vertices of $R_{2}$ opposite the axis have coordinates
$(O,E,e_{1})$ and $(O+4m,E,e_{1})$.

Let $BOX$ be that parallelepiped whose projection on planes $\mathcal{P}$ and $%
\mathcal{Q}$ is rectangles $R_{1}$ and $R_{2}$, respectively;
i.e., \[BOX =\{(x,y,z)\mid O\leqq x\leqq O+4m,E\leqq y\leqq
E+4n,e_{1}\leqq z\leqq
o_{1}\}.\] So the dimensions of $BOX$ are $4m\times 4n\times o$ where $%
o=e_{1}-o_{1}$ is odd. We assert $BOX$ is a fundamental domain for
$\widehat{G} $.

There is a tessellation of $E^{3}$ obtained by translating $BOX$
around using translations by $4m$, $4n$, and $o$ in the $x$, $y$,
and $z$ directions, respectively. One observes, from the rectangle
theorems, that the rotations in $\widehat{G}$, which generate
$\widehat{G}$, leave this tessellation invariant. Also
$\widehat{G}$ contains translations by $4m$, $4n$, and $4o$ in the
$x$, $y$, and $z$ directions, respectively. Using these
translations
and the rotations which split the faces of $BOX$ we see that any point in $%
E^{3}$ is equivalent to a point in BOX. If two points in interior
of $BOX$ are equivalent then there is a non--trivial element
$\widehat{g}$ of $\widehat{G} $ that leaves $BOX$ invariant. By
the Brouwer fixed point theorem, $\widehat{g} $ has a fixed point
in $BOX$ and therefore must be a rotation whose axis intersects
BOX. Inspecting rectangles $R_{1}$ and $R_{2}$ we see that this is
impossible. Thus $BOX$ is a fundamental domain for $\widehat{G}$.

We can conjugate $\widehat{G}$ by an element $\widehat{u}$ of
$\widehat{U}$ and obtain an equivalent subgroup of $\widehat{U}$.
This has the effect of replacing $BOX$ by $\widehat{u}(BOX)$. As
$\widehat{U}$ contains translations by 4 in the $x$, $y$, and $z$
directions we may assume without loss of generality that
$BOX=\{(x,y,z)\mid \widehat{O}-2m\leqq x\leqq
\widehat{O}+2m,\widehat{E}-2n\leqq y\leqq \widehat{E}+2n,\widehat{e}%
_{1}\leqq z\leqq \widehat{o}_{1}\}$ where $\widehat{O}=\pm 1$, $\widehat{E}%
=0 $ or 2, $\widehat{e}_{1}=0$ or 2 and $o=\widehat{o}_{1}-\widehat{e}_{1}$.
The rotations $\widehat{a}$, $\widehat{b}$, and $\widehat{c}$ of $\widehat{U}
$ are given by equations $(x,y,z)\longrightarrow (x,-y,-z+2)$, $%
(x,y,z)\longrightarrow (-x+2,,-z)$, $(x,y,z)\longrightarrow (-x,-y+2,z)$
respectively. So applying $\widehat{a}$, $\widehat{b}$, or $\widehat{c}$ if
need be we can assume $\widehat{O}=1$, $\widehat{E}=0$ and $\widehat{e}%
_{1}=0 $. But then $BOX=Box(G(m,n,o))$ which implies $\widehat{G}%
=\widehat{G}(m,n,o)$. We have shown that any finite index subgroup of $%
\widehat{U}$ generated by rotations that contains rotations with axes in
only two of the three possible directions is equivalent to a $\widehat{G}%
(m,n,o)$.

 Now suppose $\widehat{G}$ contains rotations with axes parallel to
the $x$, $y$, and $z$ directions.
For each choice of an ordered pair from the set \{$x$--axis, $y$--axis, $z$%
--axis\} to play the role of $y$--axis and $z$--axis in
Proposition 4 we get a rectangle theorem. We don't formally state
each of the six propositions but we use the results to get
tessellations of planes by rectangles in order to construct a
parallelepiped, again called $BOX$, which will turn out to be a
fundamental domain for $\widehat{G}$.

Let $\mathcal{P}$ (resp. $\mathcal{Q}$ , $\mathcal{R}$) be a plane
parallel to the $xy$ (resp. $xz$, $yz$) plane containing an axis
$a_{x}=(t$, even, odd) (resp. $a_{z}=$(even, odd,$t$),
$a_{y}=$(odd, $t$, even)) parallel to the $x$ (resp. $z$, $y$) axis. Then planes $%
\mathcal{P}$, $\mathcal{Q}$ , and $\mathcal{R}$ intersect in a point $%
X=(o_{1},o_{2},o_{3})$ with all odd coordinates. (For example, plane $%
\mathcal{P}$ contains axis $a_{x}=(t$,even, odd) and $P$ is
parallel to the $xy$ plane and so has equation $z=$ odd.) No point
with all odd coordinates belongs to an axis of rotation in
$\widehat{U}$.

Consider the tessellation of plane $\mathcal{P}$ by rectangles. Planes $%
\mathcal{P}$ and $\mathcal{Q}$ intersect in a line $\ell $ (see Figure 9)
parallel to the $x$--axis and planes $\mathcal{P}$ and $\mathcal{R}$
intersect in a line $m$ parallel to the $y$--axis. As $\mathcal{Q}$ contains
axes from $\widehat{G}$ parallel to the $z$--axis line $\ell $ contains $z$%
--axis points that are vertices of the tessellation by rectangles. We
already know that the axes in $\mathcal{P}$ parallel to the $x$--axis evenly
divide the rectangles but the line $m$ which is parallel to the $y$--axis
also evenly divides rectangles. To see this translate $\mathcal{P}$ in the $%
z $--direction to a plane $\widetilde{\mathcal{P}}$ that contains an axis
from $\widehat{G}$ that is parallel to the $y$--axis. This translation, in
the $z$--direction, takes vertices of the tessellation of $\mathcal{P}$ by
rectangles to vertices of the tessellation of $\widetilde{\mathcal{P}}$ by
rectangles, leaves plane $\mathcal{R}$ invariant and sends line $m$ to an
axis in $\widehat{G}$ parallel to the $y$--axis that evenly divides a
rectangle in $\widetilde{\mathcal{P}}$. Therefore $m$ evenly divides a
rectangle of the tessellation of $\mathcal{P}$.

The tessellations of planes $\mathcal{Q}$ and $\mathcal{R}$ by rectangles is
also displayed in Figure 9. Planes $\mathcal{Q}$ and $\mathcal{R}$ intersect
in line $n$ parallel to the $z$--axis.

\begin{figure}[ht]
\epsfig{file=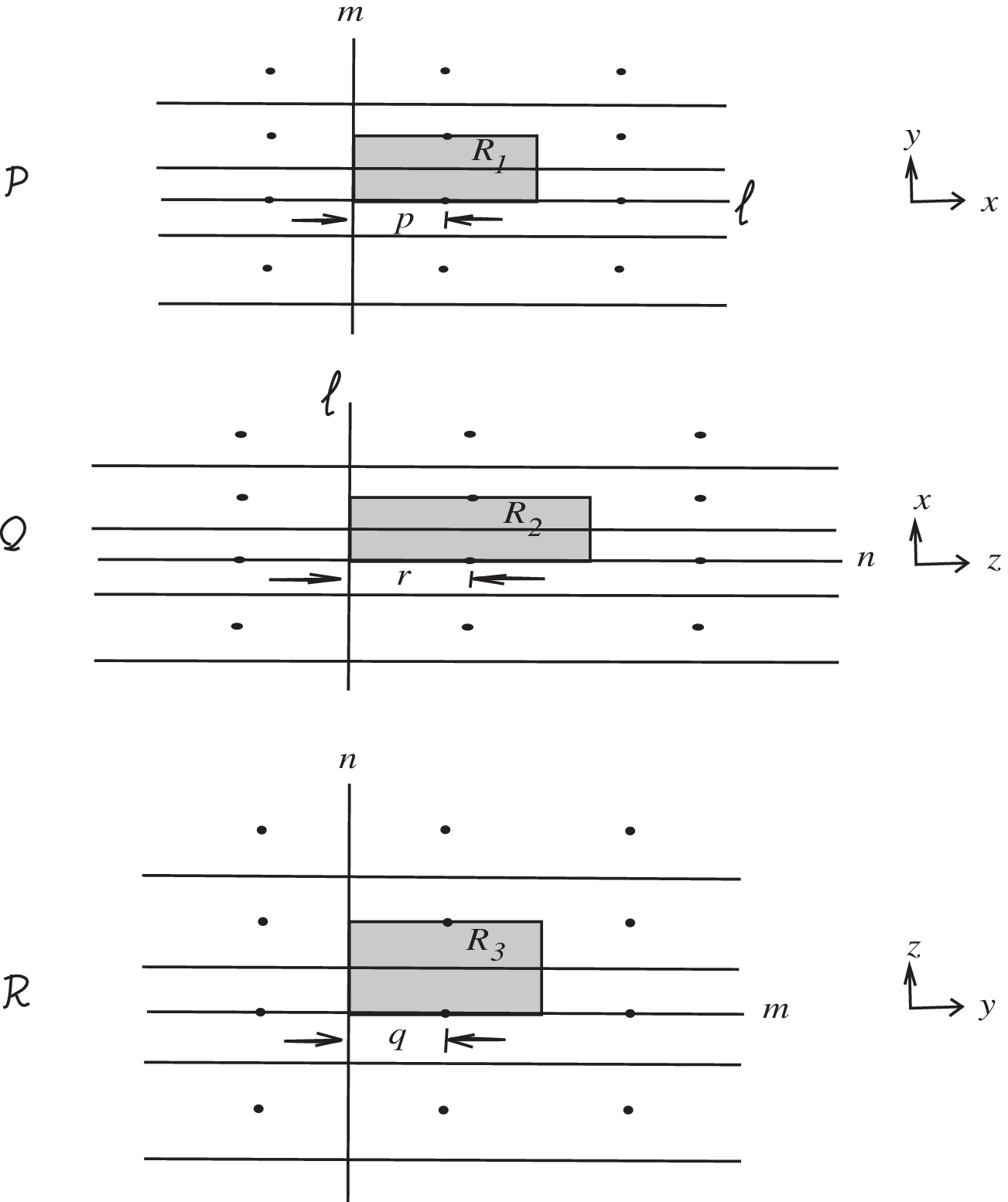,height=11.5cm}\caption{}\label{fig9}
\end{figure}

The distance from point $X=(o_{1},o_{2},o_{3})$ to the nearest axis in $%
\widehat{G}$ parallel to the $x$ (resp. $y$, $z$) axis is $q$ (resp. $r$, $p$%
) as displayed in Figure \ref{fig9}. That $p$, $q$, and $r$ are odd integers
can be seen from the parameterizations of the axes in $\widehat{U}$.

Then $BOX$ is defined to be $\{(x,y,z)\mid o_{1}\leqq x\leqq
o_{1}+2p;o_{2}\leqq y\leqq o_{2}+2q;o_{3}\leqq z\leqq o_{3}+2r\}$.
The
projections of $BOX$ on planes $\mathcal{P}$, $\mathcal{Q}$ , $%
\mathcal{R}$ are the rectangles $R_{1}$, $R_{2}$, $R_{3}$ displayed in
Figure \ref{fig9}. We assert that $BOX$ is a fundamental domain for $\widehat{G%
}$. Using translations by $4p$, $4q$, $4r$ in the $x$, $y$, and $z$
directions, which are contained in $\widehat{G}$ together with rotations,
giving rise to the vertices of the rectangles displayed in Figure \ref{fig9}%
, we see that any point in $E^{3}$ can be moved to a point in
$BOX$. There is a tessellation of $E^{3}$ obtained by translating
$BOX$ around using translations of $2p$, $2q$, and $2r$ in the
$x$, $y$, and $z$ directions. From Figure \ref{fig9} we see that
the rotations in $\widehat{G}$ preserve this tessellation so that
$\widehat{G}$ preserves this tessellation. If two points in the
interior of $BOX$ are equivalent, say $g(m_{1})=m_{2}$, then $g$
preserves $BOX$, has a fixed point by the Brouwer Theorem and so
must be a rotation in $\widehat{G}$ which is impossible.

Therefore $BOX$ is a fundamental domain for $\widehat{G}$. As
before
we can conjugate by an element of $\widehat{u}$ of $\widehat{U}$ or $%
\widehat{S}$ which has the effect of replacing $\widehat{G}$ by an
equivalent group and $BOX$ by the new fundamental domain
$\widehat{u}$(BOX).
 Since the center of $BOX$ has all even coordinates $%
(o_{1}+p,o_{2}+q,o_{3}+r)$ we can find an element $\widehat{u}$
which is a product of translations by 4 in the $x$, $y$, and $z$
directions and
rotations $\widehat{a}$ and/or $\widehat{b}$ and/or $\widehat{c}$ such that $%
\widehat{u}$(BOX) is centered at the origin. Finally we can use
the 120 degree rotation $D$ in $\widehat{S}$ to cyclically permute
$p$, $q$, $r$ so that $p\leqq \max \{q,r\}$ and in the case where
$p$, $q$, and $r$ are not all different $p\leqq q\leqq r$. Thus
the group to which $\widehat{G}$ is equivalent is
$\widehat{H}(p,q,r)$ which has $\widehat{u}$(BOX) as its
fundamental domain. We summarize all this as a theorem.

\begin{theorem}
\label{t5} 1. Let $\widehat{G}$ be an even index subgroup of $\widehat{U}$
generated by rotations. Then $\widehat{G}$ is equivalent to a unique group
in the family $\widehat{G}(m,n,o)$ and $\big[\widehat{U}:\widehat{G}\big]%
=2mno$. $\widehat{G}$ contains axes of rotation in two of the three
directions $x$, $y$, and $z$. The integers $2m$ (resp. $2n$) represents the
distance between adjacent axes of $\widehat{G}(m,n,o)$ that lie in a plane
parallel to the $xy$ plane and are parallel to the $y$--axis (resp. $x$
axis). The odd integer $o$ represents this distance between axes of $%
\widehat{G}$ that are not parallel but are as close as possible.

2. Let $\widehat{G}$ be an odd index subgroup of $\widehat{U}$ generated by
rotations. Then $\widehat{G}$ is equivalent to a unique group in the family $%
\widehat{H}(p,q,r)$ and $\big[\widehat{U}:\widehat{G}\big]=pqr$. Group $%
\widehat{G}$ contains rotations with axes parallel to each of the three
possible directions $x$, $y$, and $z$.
\end{theorem}

For each pair of directions $x$ and $y$, $x$ and $z$, $y$ and $z$ there is a
distance between a pair of axes in these directions that are not parallel
but are as close as possible giving rise to a triple of integers. This
triple of integers is $p$, $q$, and $r$, not necessarily in that order.

In the next section, we begin the study of finite index subgroups of $U$
that are generated by rotations.

\section{Finite index subgroups of $U$ generated by rotations}

\begin{proposition}
\label{p6} Let $\widehat{G}$ be a finite index subgroup of $\widehat{U}$
generated by rotations and let $G=\varphi ^{-1}(\widehat{G})$ be the full preimage of $%
G$ under the homomorphism $\varphi :U\rightarrow \widehat{U}$
defined in Section 2. Then $G$ is generated by rotations.
\end{proposition}

\begin{proof}
The homomorphism $\varphi :G\rightarrow \widehat{G}$ defined in
Section 2 is surjective, and sends $90^{\circ}$ rotations in $U$
to $180^{\circ}$ rotations in $\widehat{U}$. By the classification
of the $\widehat{G}$ in Section 4, $\widehat{G}$ is generated by 3
or 4 rotations. Let $S$ be a set
of $90^{\circ }$ rotations in $U$ that is sent to a set of generators for $%
\widehat{G}$ and let $G_{1}$ be the subgroup of $U$ generated by $S$. Then $%
\varphi ^{-1}(\widehat{G})=G_{1}K$. Since $K$ is generated by rotations
(Theorem 2), so is $G_{1}K$.
\end{proof}

The main theorem now follows easily from Proposition \ref{p6}.

\begin{theorem}
\label{t7} Given any positive integer $n$ there is a subgroup $G$ of $U$ of
index $n$ that is generated by rotations.
\end{theorem}

\begin{proof}
Let $\widehat{G}$ be a subgroup of $\widehat{U}$ generated by rotations of
index $n$ in $\widehat{U}$, which exists by the classification of such
subgroups of Section 4. And let $G=\varphi ^{-1}(\widehat{G})$. Then $G$ is
generated by rotations by Proposition 6 and $[U:G]=[\widehat{U}:\widehat{G}%
]=n$.
\end{proof}

Any axis of rotation $\ell$ in $U$ is the image of the axis of rotation of
one of the generators $a$, $b$, $c$ of $U$ under the action of an element $u$
of $U$. This follows from the fact that $D_0$, a dodecahedral fundamental
domain of $U$, intersects six axes of rotation in $U$, those of $a$, $b$, $c$%
, $c^{-1}ac$, $a^{-1}ba$, and $b^{-1}cb$, and if $D$ is any dodecahedron of
the tessellation of $H^{3}$ intersecting $\ell$ there is an element $u_1$ of
$U$ such that $u_1(D_0) = D$. Then $u = u_1x^{-1}$ where $x$ is one of $a$, $%
b$, $c$. Letting $U$ act on the axes of rotation, we get exactly three
orbits. (At most three by the argument above and at least three because $%
\varphi : U \longrightarrow \widehat U$ preserves orbits and there are three
orbits in $\widehat U$, those parallel to the $x$, $y$ and $z$ axes.) Thus
there are nine conjugacy classes of rotations in $U$ represented by $a$, $a^2
$, $a^3$, $b$, $b^2$, $b^3$, $c$, $c^2$, and $c^3$. (This can also be seen
by computing $U/[U, U] \cong Z_4 \oplus Z_4 \oplus Z_4$ from the
presentation of $U$ in Section 2. For example $a$ is sent to $(1, 0, 0)$,
etc.) Similarly there are three conjugacy classes of rotations in $\widehat U
$ represented by $\widehat a$, $\widehat b$, and $\widehat c$.

\begin{theorem}
\label{t8} Let $G$ be a subgroup of $U$ of odd index and generated by
rotations. Then $G$ contains a member of each of the nine conjugacy classes
of rotations in $U$.
\end{theorem}

\begin{proof}
Let $\widehat{G}=\varphi (G)$ where $\varphi :U\rightharpoonup \widehat{U}$
is the homomorphism of Section 2 and $K=\ker \varphi $. Then $G\subset
GK\subset U$ so that $[U:G]=[U:GK][GK:G]$. But $\varphi :U\longrightarrow
\widehat{U}$ induces $\varphi :GK\longrightarrow \widehat{G}$ so that $%
[U:GK]=[\widehat{U}:\widehat{G}]$ and $[U:G]=[\widehat{U}:\widehat{G}]\cdot
\lbrack GK:G]$. Since $[U:G]$ is odd it follows that $[\widehat{U}:\widehat{G%
}]$ is odd and thus $\widehat{G}$ contains a member of each of the three
conjugacy classes of rotations in $\widehat{U}$ from the classification in
Section 4.

We shall show that $G$ contains a member of the conjugacy class of $c$. Let $%
\widehat{c}_{1}$ be a rotation in $\widehat{G}$ with axis parallel to the $z$%
--axis. Suppose $\varphi (g)=\widehat{c}_{1}$. Then $g$ is a product of
rotations, $g=\prod\limits_{i=1}^{n}r_{i}$, as $G$ is generated by
rotations. If $\{r_{1},\dots ,r_{n}\}$ contains a rotation conjugate to $c$
or $c^{3}$ we are done. Suppose this is not the case. Then $\widehat{c}%
_{1}=\prod\limits_{i=1}^{n}\widehat{r}_{i}$ where $\widehat{r}_{i}$ is
either the identity or a rotation about an axis parallel to the $x$ or $y$
axes. Each $\widehat{r}_{i}$ belongs to the group $\widehat{G}(1,1,1)$
defined in Section 4 as that group contains every rotation in $\widehat{U}$
about an axis parallel to the $x$ or $y$ axis. Thus $\widehat{c}_{1}\in
\widehat{G}(1,1,1)$ but this is impossible as $\widehat{G}(1,1,1)$ contains
no rotations with axis parallel to the $z$--axis. Therefore $G$ contains a
member of the conjugacy class of $c$.

The two conjugates $D\widehat{G}(1,1,1)D^{-1}$ and $D^{2}\widehat{G}%
(1,1,1)D^{-2}$, where $D$ is $120^{\circ }$ rotation about axis $(t,t,t)$
introduced in Section 4, contain all rotations parallel to the $x$ and $z$
axis and no rotation parallel to the $y$ axis or all rotations parallel to
the $y$ and $z$ axes and no rotations parallel to the $x$ axis. We can show
that $G$ contains rotations in the conjugacy class of $a$ and $b$ by
duplicating the argument for $c$ by replacing $\widehat{G}(1,1,1)$ by $D%
\widehat{G}(1,1,1)D^{-1}$ or $D^{2}\widehat{G}(1,1,1)D^{-2}$.
\end{proof}

If $G$ is a finite index subgroup of $U$ that is generated by rotations it
is clear that information about the precise placement of $\widehat G$ in the
classification of Section 4 implies much about group $G$ itself. There are
other theorems analogous to Theorems 7 and 8, but clumsier to state or prove
that we could present. We refrain from doing so, so as not to lengthen this
paper.

We close by posing a question. If $G$ is a subgroup of $U$ of index $n$,
either generated by rotations or not, it is clear that $G$ has a fundamental
domain that is a union of $n$ of the dodecahedra in the tessellation
associated to $U$. Does $G$ have a fundamental domain that is convex and
also the union of $n$ dodecahedra?

\end{document}